\documentclass[12pt]{amsart}

\def \mindot{{\dot \circ}}
\def \mincross{{\dot \times}}
\def \fg{{\mathfrak g}}
\def \fh{{\mathfrak h}}
\def \fk{{\mathfrak k}}

\def \ft{{\mathfrak t}}

\def \u{{\mathfrak u}}
\def \su{{\mathfrak s}{\mathfrak u}}
\def \gl{{\mathfrak g}{\mathfrak l}}
\def \C{{\mathbb C}}
\def \T{{\mathbb T}}
\def \R{{\mathbb R}}
\def \P{{\mathbb P}}
\def \H{{\mathcal H}}

\def \W{{\mathcal W}}
\def \OO{{\mathcal O}}

\newtheorem{theorem}{Theorem}[section]

\newtheorem{proposition}[theorem]{Proposition}

\numberwithin{equation}{section}

\begin{document}

\baselineskip=16pt

\title[Polygons in Minkowski space]{Polygons in Minkowski space and Gelfand-Tsetlin for pseudounitary 
groups}

\author[P. Foth]{Philip Foth}

\address{Department of Mathematics, University of Arizona, Tucson, AZ 85721-0089}

\email{foth@math.arizona.edu}



\date{April 14, 2007}

\begin{abstract}
 
We study the symplectic geometry of the moduli spaces of polygons in the Minkowski 3-space.
These spaces naturally carry completely integrable systems with periodic flows.
We extend the Gelfand-Tsetlin method to pseudo-unitary groups and show that the action
variables are given by the Minkowski lengths of non-intersecting diagonals.  
\end{abstract}

\maketitle

\section{Introduction} The geometry of moduli spaces of polygon in the Eucledian 3-space 
has been studied by many authors, notably by Klyachko \cite{Klyachko}, Kapovich and Millson \cite{KM}, 
Haussman and Knutson \cite{HK}, among others, and many interesting results in 
symplectic geometry have been obtained. In this paper we study polygons in the Minkowski 
3-space and obtain a variety of results, similar in spirit, but which are, on the other hand, 
considerably different, and the differences are illuminated by the use of pseudounitary groups
${\rm U}(p,q)$ and their coadjoint orbits. 

Polygons in the Minkowski 3-space were briefly considered by Millson 
in \cite{Millson}. Original application of the Gelfand-Tsetlin method to integrable systems on coadjoint
orbits is due to Guillemin and Sternberg \cite{GS}.

\section{Polygons in Minkowski space}

A surface $\H_R$ in $\R^3$ defined by the equation $t^2-x^2-y^2=R^2$ is called a {\it pseudosphere}, 
since the Minkowski metric of signature (2,1) restricts to the constant curvature Riemannian 
metric on it. Alternatively, we can think of a pseudosphere as a set of points equidistant from
the origin in $\R^3$ with respect to the Minkowski metric. The connected component $\H^+_R$ corresponding
to $t >0$ will be called a future pseudosphere, and $\H^-_R$ corresponding to $t<0$ a past pseudosphere respectively.
Note that the group ${\rm SU}(1,1)$ acts transitively on each connected component, since we can think
of $\R^3$ as $\su(1,1)^*$, where these connected components can be
thought of as elliptic coadjoint orbits. Therefore
each has a natural symplectic structure, invariant under the action of the group. The metric 
is invariant as well, since ${\rm SU}(1,1)$ acts by isometries. Each connected component is 
also a K\"ahler manifold, since it is isomorphic to the hyperbolic plane ${\rm SU}(1,1)/{\rm U}(1)$.  

The purpose of this section is to study the geometry of the symplectic quotients of the product of 
several future and past pseudospheres with respect to the diagonal ${\rm SU}(1,1)$-action. These
spaces have a natural interpretation as polygon spaces in Minkowski 3-space. Let us start by fixing 
notation. Let ${\bf r}=(r_1, ..., r_n)$ be an $n$-tuple of positive real numbers, and let us fix two 
positive integers $p\ge q$ such that $p+q =n$. In the language of polygons, this will mean that we
have the first $p$ sides in the future timelike cone and the last $q$ in the past. The Minkowski
length of the $i$-th side is equal to $r_i$ and the space of closed polygons, i.e. those where the
sum of the first $p$ sides in the future timelike cone equals the negative of the sum of 
the last $q$ sides in the past timelike cone, is identified with the zero level set of the moment map:
$$
\mu:\ \OO_{1}\times \cdots \times \OO_n \to \su(1,1)^*\ .
$$
Here $\OO_i\simeq \H^+_{r_i}$ is a future pseudosphere of radius $r_i$ if $1\le i\le p$ and past 
$\OO_i\simeq \H^-_{r_i}$ if 
$p+1\le i\le n$ with its coadjoint orbit sympletic structure. Note that the triangle inequalities in 
the future (or past) timelike cone are reversed from the usual ones. If $v_1$ and $v_2$ are two 
equally directed timelike vectors, then $||v_1 +v_2||\ge ||v_1|| + ||v_2||$. For convenience, let us 
fix the perimeter of the polygon to be equal to $2$. This means that $\sum_{i=1}^n r_i=2$. Note that 
each $\OO_i$ can itself be naturally interpreted as the symplectic quotient of $\C^2$ with complex 
coordinates $(z, w)$, and symplectic form 
$\displaystyle{\frac{\sqrt{-1}}{2}(dz\wedge d{\bar z} - dw\wedge d{\bar w})}$ 
for $1\le i\le p$ and negative of that
for $p+1\le i\le n$, at the level $r_i$, with respect to the diagonal circle action. 

Let us denote by $M_{\bf r}$ the quotient $ \mu^{-1}(0)/{\rm SU}(1,1)$. This is a quotient of a 
non-compact space, in general, by the action of a non-compact Lie group. Therefore, questions of its
topology and geometry require careful consideration. However, notice that for a generic choice of
${\bf r}$, i.e. such that $M_{\bf r}$ is non-empty and $r_1+\cdots +r_p\ne r_{p+1}+\cdots +r_n$,
every point in the moduli space represents a polygon with a trivial stabilizer. This means that the
group ${\rm PSU}(1,1)$ acts freely and properly on $\mu^{-1}(0)$. An easy computation shows that with our
assumptions, $0$ is a regular value for the map $\mu$, and therefore \cite[Theorem 1.11.4]{DK}, 
the quotient space $M_{\bf r}$ 
has the structure of a smooth manifold of dimension $2n-6$. In a later section we will see a different 
approach to $M_{\bf r}$ through a fixed point set of an involution in a hyper-K\"ahler manifold.

One of the powerful tools in dealing with polygons in a compact setting proved to be \cite{HK} reduction 
in stages, or symplectic Gelfand-MacPherson correspondence. In fact, an appropriate modification of this
method proves to be useful for our purposes as well. 

Let us consider the space $\C^{2n}$ with complex coordinates $(z_1, ...,z_n,w_1, ..., w_n)$ and symplectic 
form $\Omega$ given by 
$$
2\sqrt{-1}\Omega = \sum_{i=1}^p (dz_i\wedge d{\bar z}_i - dw_i\wedge d{\bar w}_i) -
\sum_{j=p+1}^n (dz_j\wedge d{\bar z}_j - dw_j\wedge d{\bar w}_j).
$$
We introduce two elements ${\bf z}= (z_1, ..., z_n)^T$ and ${\bf w}=(w_1, ..., w_n)^T$ of $\C^n$ and comprise 
an $n\times 2$ matrix $M=({\bf z}\ {\bf w})$, representing an element of $\C^{2n}$. There is a natural 
left action of ${\rm U}(p,q)$ on $\C^{2n}$ given by left multiplication $M\mapsto AM$, where 
$A\in{\rm U}(p,q)$ and $M$ as before, and similarly a natural right action of ${\rm U}(1,1)$. Both are 
hamiltonian actions with respective moment maps 

$$
\eta: \C^{2n}\to \u(p,q)^*, \ \ \eta(M)=MJ_{1,1}M^*J_{p,q} \ ,
$$
 
where $M^* = {\bar M}^T$, as usual, and

$$
\nu: \C^{2n}\to \u(1,1)^*, \ \ \nu(M)=J_{1,1}M^*J_{p,q}M \ ,
$$
where $J_{1,1}={\rm diag}(1,-1)$ and 
$J_{p,q}={\rm diag}(\underbrace{1, ..., 1}_p, \underbrace{-1,  ..., -1}_q)$.

To make it more explicit, we note that 
$$
\eta(M) = \left( \begin{array}{cc}
||{\bf z}||^2 & \langle {\bf w}, {\bf z} \rangle \\
- \langle {\bf z}, {\bf w} \rangle & - ||{\bf w}||^2
\end{array}
\right) \ ,
$$
where the norm and the pairing come from the standard pseudohermitian structure on $\C^n$ 
of signature $(p,q)$. 

We notice that the left action of the diagonal $\T^n\subset {\rm U}(p,q)$ commutes with the 
right action of ${\rm U}(1,1)$.  
Now, shifting for convenience by the identity matrix (i.e. the central matrix with the trace
equal to the perimeter), we can look at the level set of $\nu$ corresponding to the identity
$2\times 2$ matrix $I_2$ in $\u(1,1)^*$, i.e. the orthonormal pairs of vectors 
$({\bf z}, {\bf w})$ in $\C^n$ such that ${\bf z}$ is timelike and ${\bf w}$ is spacelike. 
The quotient of this level set by the aforementioned action of ${\rm U}(1,1)$ is naturally 
isomorphic to the semisimple symmetric space 
$$
X_{p,q,1,1} \simeq {\rm U}(p,q)/{\rm U}(1,1)\times {\rm U}(p-1, q-1)
$$
(recall that we are working under the assumption that $p\ge q\ge 1$). This space is a pseudo-hermitian 
symmetric space; it has invariant complex and compatible symplectic structures. 

The residual hamiltonian action of $\T^n$ on $X_{p,q,1,1}$ has as a moment map $\eta_\T$ the projection of 
$\eta(M)$ onto the diagonal, i.e. 
$$
\eta_\T(M) = (|z_1|^2-|w_1|^2, ..., |z_p|^2-|w_p|^2, |w_{p+1}|^2 - |z_{p+1}|^2, ..., |w_n|^2-|z_n|^2). 
$$ 
Naturally, the quotient of the level set of $\eta_\T$ corresponding to $(r_1, ..., r_n)$ is the moduli
space of polygons in question $M_{\bf r}$. 

This correspondence between the two symplectic quotients helps to understand the following important
feature of the space $M_{\bf r}$:

\begin{proposition}
If $q=1$, the space $M_{\bf r}$ is compact, and if $q>1$, the space $M_{\bf r}$ is not compact. 
\end{proposition}
\noindent{\it Proof.} The coadjoint orbit of $\eta(M)$ is elliptic, and passes through 
$$\Lambda = {\rm diag}(1, \underbrace{0, ..., 0}_{p-1}, 1, \underbrace{0, ...0}_{q-1}).$$ We will show that 
for $q>1$ there is no positive adapted system of roots, in terminology of \cite[Definition VII.2.6]{Neeb}, for 
which $\Lambda$ is admissible. Then [{\it loc.cit.}, Theorem VIII.1.8] would immediately imply that the
map $\eta_\T: X_{p,q,1,1}\to \ft^*$ is not proper. Given that ${\bf r}$ is chosen generic in the image, 
we will be able to conclude the statement. However, one characterization of the positive adapted root system,
[{\it loc.cit.}, Proposition VII.2.12] valid for quasihermitian Lie algebras, implies that the condition of 
$\Delta^+$ being adapted is equivalent to the system $\Delta^+_n$ of positive non-compact roots being 
invariant under the baby Weyl group. Now it is easy to see that the latter is possible if and only if 
$q=1$, in which case $\Delta^+$ should be taken the negative of the standard subset of positive roots.
{\bf Q.E.D.} \medskip

Another, more visual way of seeing that $M_{\bf r}$ is only compact when $q=1$, can be found using polygons. 
First, let us explain compactness for $q=1$. The last side of the polygon, ${\bf e}_n$, of Minkowski length
$r_n$ can be represented, after applying the action by an element of ${\rm SU}(1,1)$ by a vector in $\R^3$
with coordinates $(0,0,-r_n)$. Therefore, the $(n-1)$ future timelike sides of the polygon should add 
up to $(0,0,r_n)$, since the only degree of symmetry left is the circle rotation around the $t$-axis. 
Clearly, this space is bounded and closed, therefore compact. On the contrary, when $q>1$, the space $M_{\bf r}$
is not compact. Let us explain this in the simplest example $p=q=2$ and $r_1=r_2=r_3=r_4=1/2$. Note that
again, using the action of ${\rm SU}(1,1)$, we can assume that ${\bf e}_4=(0,0,-1/2)$, and the only degree 
of symmetry left is again the rotation about the $t$-axis. We will produce an explicit sequence 
of points in $M_{\bf r}$ with no limit. Let $x_n$ be the closed polygon corresponding to 
${\bf e}_1=-{\bf e}_4=(0,0,1/2)$, ${\bf e}_2=(n, 0, \sqrt{n^2+1/4})$, and ${\bf e}_3=-{\bf e}_2$. 
Clearly, the sequence $(x_n)$ has no limit points in $M_{\bf r}$ and thus $M_{\bf r}$ is not compact.

Let us denote by $d_i$ the length of $i$-th diagonal, i.e. the diagonal connecting the first and the 
$(i+1)$-st vertex. By our convention, $d_1 = r_1$, $d_{n-1}=r_n$, and $d_n=0$. In the next section 
we will show that similarly to the compact situation,  the lengths of the $(n-3)$ varying diagonals 
$d_2$, ..., $d_{n-2}$ define a completely integrable system on $M_{\bf r}$, and are action variables 
for $(n-3)$ periodic flows. 

\section{Symplectic structure on the moduli space}

In this section we will spell out the elementary definition of the symplectic structure on the space 
$M_{\bf r}$, in quite a similar way to \cite[Section 3]{KM}. First of all, let us define the following
two operations on $\R^3$, with coordinate functions $(x,y,t)$. For two vectors 
${\bf v}_1=(x_1, y_1, t_1)^T$ and ${\bf v}_2 = (x_2, y_2, t_2)^T$ we define the Minkowski cross product 
$\mincross$ and the Minkowski dot product $\mindot$ as follows:
$$
{\bf v}_1\mincross{\bf v}_2 = {\rm det}\left(
\begin{array}{ccc}
-{\bf i} & -{\bf j} & {\bf k} \\ x_1 & y_1 & t_1 \\ x_2 & y_2 & t_2  
\end{array} 
\right)\ , \ \ {\rm and} \ \ 
{\bf v}_1\mindot{\bf v}_2 = -x_1x_2 - y_1y_2 +t_1 t_2 \ \ ,
$$
where ${\bf i}$, ${\bf j}$, and ${\bf k}$ are the usual unit vectors in the positive directions
of the $x$-, $y$-, and $t$-axes respectively. Note that $\mindot$ is non-degenerate and 
positive definite in the timelike cone. 

These operations satisfy the usual properties of the dot and cross products in $\R^3$:

$$
\begin{array}{c}
{\bf a}\mincross{\bf b} = - {\bf b}\mincross{\bf a} \\
({\bf a}\mincross{\bf b})\mincross{\bf c}+({\bf b}\mincross{\bf c})\mincross{\bf a}+({\bf c}\mincross{\bf a})\mincross{\bf b}=0 \\
{\bf a}\mincross({\bf b}\mincross{\bf c})={\bf b}({\bf a}\mindot{\bf c})-{\bf c}({\bf a}\mindot{\bf b}) \\
{\bf a}\mindot ({\bf b}\mincross{\bf c})={\rm det}({\bf a}\ {\bf b}\ {\bf c}) \\
({\bf a}\mincross {\bf b})\mincross ({\bf c}\mincross{\bf d}) =
{\rm det}({\bf a}\ {\bf b}\ {\bf d}){\bf c} - {\rm det}({\bf a}\ {\bf b}\ {\bf c}){\bf d} 
\end{array}
$$
The first two properties show that $(\R^3, \mincross)$ is a Lie algebra, in fact isomorphic to $\su(1,1)$
under the following map:
$$
\left(
\begin{array}{c} x \\ y \\ t \end{array} 
\right) \mapsto \frac{1}{2}\left(
\begin{array}{cc}
-\sqrt{-1}\cdot t & x+\sqrt{-1}\cdot y \\ x-\sqrt{-1}\cdot y & \sqrt{-1}\cdot t
\end{array}
\right) \ .
$$
Under this identification, $\mindot$ corresponds to $-2{\rm Tr}(AB)$. 

Now the description of the symplectic two-form $\omega$ on the hyperboloid $\H_R$ given by the 
equation $t^2-x^2-y^2=R^2$ is given by 
$$
\omega_{\bf u}({\bf v_1}, {\bf v}_2) = \frac{1}{R^2}{\bf u}\mindot({\bf v}_1\mincross{\bf v}_2) \ ,
$$
where ${\bf u}$ is a point on the hyperboloid, and ${\bf v}_1$ and ${\bf v}_2$ are 
elements of ${\rm T}_{\bf u}\H_R$. Here we think of ${\rm T}_{\bf u}\H_R$ as the linear subspace
of $\R^3$ orthogonal to ${\bf u}$ with respect to $\mindot$. 
Similarly to \cite[Lemma 3.1]{KM} we see that the map 
$$
\begin{array}{c}
\H_{r_1}\times \cdots \times \H_{r_n} \to \R^3 \\
({\bf u}_1, ..., {\bf u}_n)\mapsto {\bf u}_1+\cdots +{\bf u}_n
\end{array}
$$
is the moment map with respect to the diagonal ${\rm SU}(1,1)$-action and the product 
symplectic structure. 

Now we will describe the hamiltonian flow $\phi_i(t)$ on the space $M_{\bf r}$ corresponding 
to the hamiltonian function $d_i$ - the Minkowski length of the $i$-th diagonal of the polygon,
connecting the first and the $(i+1)$-st vertices. Note that 
$$
d_i^2 = {\bf d}_i\mindot {\bf d}_i
=({\bf u}_1+\cdots +{\bf u_{i+1}})\mindot  ({\bf u}_1+\cdots +{\bf u_{i+1}})
$$
is a positive real number, since the vector in parenthesis is in the (future) timelike cone, 
by our assumptions, so we take $d_i$ real positive as well. Note that if we place the first
vertex at the origin of $\R^3$ and use the action of ${\rm SU}(1,1)$ to move the $(i+1)$-st
vertex to a position on the $t$-axis, then the corresponding bending flow is easy to describe
as rotation of vertices numbered $2$, ..., $i$ about the $t$-axis with a constant angular 
speed, since the Hamiltonian vector field in this case is given by 

\begin{equation}
({\bf d}_i\mincross {\bf u}_1, ..., 
{\bf d}_i\mincross {\bf u}_{i}, 0, ..., 0) = d_i(y_1{\bf i}-x_1{\bf j}, ...,
y_i{\bf i}-x_i{\bf j}, 0, ...,0)\ . 
\label{eq257}
\end{equation}

The general statement follows from the equivariance with respect to the ${\rm SU}(1,1)$-action.
This shows that the flows are indeed periodic with periods equal to $2\pi/d_i$. 

Next, we wish to describe the angle variables, which, however clear are from the preceding 
description, can be further illuminated by the formula analogous to 
\cite[Equation 7.1.1]{FM} for the compact case: 
$$
\cos \phi_i = \frac{({\bf d}_i\mincross {\bf u}_i)\mindot({\bf d}_i\mincross {\bf u}_{i+1})}
{||{\bf d}_i\mincross {\bf u}_i||\cdot ||{\bf d}_i\mincross {\bf u}_{i+1}||}
$$
Let us explain why this formula is true. We can assume, as before, that the diagonal 
${\bf d}_i$ is aligned with the positive direction of the $t$-axis, i.e. ${\bf d}_i=d_i{\bf k}$.
According to formula (\ref{eq257}), both Minkowski cross products in the numerator are
in the $xy$-plane, and thus the Minkowski dot product is just the negative of the 
usual Eucledian dot product. Now, the denominator has Minkowski norms of two vectors in
the $xy$-plane, each of which equals $\sqrt{-1}$ times the Eucledian norm. Therefore, 
the expression yields the cosine of the oriented dihedral angle between the two planes, which
is the $i$-th angle variable. Obviously, this formula holds in general as well, since it is 
invariant under the action of ${\rm SU}(1,1)$.

\section{Gelfand-Tsetlin system for ${\rm U}(p,q)$}

Let $p$ and $q$ be positive integers, $p\ge q$, $n=p+q$. For the Lie algebra $\fg=\u(p,q)$
we use the form ${\rm Tr}(AB)$ to identify its dual space $\u(p,q)^*$ with $\sqrt{-1}\cdot\u(p,q)$, 
i.e. the space of $n\times n$ matrices $A$ such that $JA$ is Hermitian symmetric, where 
$J = {\rm diag}(\underbrace{1,1,...,1}_p, \underbrace{-1, -1, ..., -1}_q)$. In the block form, 
$$
A=\left(
\begin{array}{cc}
H_p & B \\ 
{} & {} \\
-{\bar B^T} & H_q
\end{array}
 \right)\ \ ,
$$
where $H_p$ and $H_q$ are $p\times p$ and $q\times q$ Hermitian symmetric matrices 
respectively and $B$ is a complex $p\times q$ matrix. The complexification of $\u(p,q)$
is, as usual, $\fg_{\C} = \gl(n, \C)$, and let us denote by $\fh$ the diagonal Cartan subalgebra, 
the complexification of the compact Cartan $\ft$ in $\u(p,q)$. Let $\Delta$ be the root system 
with respect to $(\fg, \fh)$, let $\Delta^+$ be the standard subset of positive roots of the 
form $e_{ij}$, $1\le i < j \le n$. We note that the positive roots $e_{ij}$ are compact if 
$1\le i,j \le p$ or $p < i,j \le n$ and non-compact if $i\le p < j$. Denote by $\Delta^+_{n}$ 
the subset of positive non-compact roots, by $\W$ the Weyl group and by $\W_c$ the baby Weyl group
corresponding to the pair $(\fk, \fh)$, where $\fk = \u(p)\times \u(q)$ is the standard maximal 
compact subalgebra in $\fg$. 

Let us now choose an $n$-tuple of real numbers $$\Lambda = (\lambda_1, ...\lambda_p, \mu_1, ..., \mu_q)$$
satisfying the following conditions: $\lambda_i \le \lambda_j$ for $i < j$ and 
$\mu_i \le \mu_j$ also for $i < j$. Besides, we will require that $\lambda_1 > \mu_q$. We will refer
to this property as $\Lambda$ being {\it admissible}. 
The condition of admissibility 
is sufficient to guarantee that the Gelfand-Tsetlin variables for the elliptic coadjoint orbit 
$\OO_\Lambda:= G.\Lambda$ are all real (here we consider $\Lambda$ as a diagonal matrix representing 
an element of $\ft^*\subset \fg^*$). Later on, we will consider a particular case, relevant to polygons, 
where this condition does not hold, but where the Gelfand-Tsetlin variables are still real. 
By \cite[Theorem 5.17]{HNP}, the projection of $\OO_\Lambda$ onto $\ft^*$ is the sum of the convex polyhedron 
${\rm conv}(\W_c.\Lambda)$ and the convex polyhedral cone defined by the non-compact positive roots. 

\medskip
\noindent{\it Remark.} 
It is an interesting separate question to find a necessary condition, which would force the eigenvalues 
of a truncated matrix to be real as well. A mere requirement that $\Lambda$ is elliptic does not
guarantee this, as we found using a computer algebra software, already in the case of ${\rm SU}(2,2)$.
\medskip

Let $\fg_{n-1}$ be the subalgebra of $\fg$ corresponding to the left upper principal submatrix of size 
$(n-1)\times (n-1)$. The algebra $\fg_{n-1}$ is isomorphic to $\u(p, q-1)$. Denote by $p_n$ the projection
$\fg^*\to \fg_{n-1}^*$. The image $p_n(\OO_\Lambda)$ is the union of certain coadjoint orbits of $U(p,q-1)$ 
in $\fg_{n-1}^*$, which we will describe next. For convenience, we denote  
$$
{\bf x}^\dagger = (J{\bar{\bf x}})^T \ .
$$

\begin{proposition} If a coadjoint ${\rm U}(p,q-1)$-orbit $\OO$ 
is in the image $p_n(\OO_\Lambda)$, then 
it is admissible. If we arrange the real eigenvalues of $\OO$ in the non-decreasing order 
$$\mu^{n-1}_1\le \mu^{n-1}_2 \le .. \le \mu^{n-1}_{q-1} \le \lambda^{n-1}_1 
\le \lambda^{n-1}_2 \le ... \le \lambda^{n-1}_p ,$$
then the following interlacing conditions hold: 
for $1\le i \le q-1$, $\mu_i \le \mu^{n-1}_i \le \mu_{i+1}$, also for $1\le j \le p-1$, 
$\lambda_j\le \lambda^{n-1}_j\le \lambda_{j+1}$ and $\lambda^{n-1}_p \ge \lambda_p$.   
\end{proposition}

\noindent{\it Proof.} The admissibility directly follows from \cite[Theorem 5.17]{HNP}.
To show interlacing, one should adapt the Courant-Fischer Theorem
\cite{Horn} separately to the timelike cone and to the spacelike cone. For example, if we let 
${\bf v}_i$ be the eigenvector in the timelike cone of $\C^n$ for the eigenvalue 
$\lambda_i$ and ${\bf w}_j$ in the spacelike cone for $\mu_j$, then 
for $$
{\bf x}=\alpha_1{\bf v}_1 +\cdots + \alpha_p {\bf v}_p 
+\beta_1{\bf w}_1 + \cdots \beta_q{\bf w}_q
$$
the Rayleigh 
quotient modifies to 
$$
R_A({\bf x}) = \frac{{\bf x}^\dagger A{\bf x}}{{\bf x}^\dagger{\bf x}}
= \frac{\sum_{i=1}^p |\alpha_i|^2\lambda_i - \sum_{j=1}^q|\beta_j|^2\mu_j}
{\sum_{i=1}^p |\alpha_i|^2 - \sum_{j=1}^q|\beta_j|^2} \ ,
$$
from which one concludes that in the timelike and spacelike cones $\C^n_+$ and $\C^n_-$ we 
respectively have  
$$
\lambda_1 = \min_{{\bf x}\in\C^n_+} R_A({\bf x})\ \ {\rm and} \ \ 
\mu_q = \max_{{\bf x}\in\C^n_-} R_A({\bf x}) .
$$
Continuing with the standard minmax arguments for timelike and spacelike cones, one gets the 
full set of interlacing conditions. In order to show that each interlacing pattern can be 
obtained this way, one can just suitably adapt the arguments of \cite[Theorem 4.3.10]{Horn}
to the pseudo-Hermitian case at hand. {\bf Q.E.D.} \medskip  
 
\noindent{\it Remark.} The Gelfand-Tsetlin patterns for the unitary representations of 
${\rm U}(p,q)$ with highest weights were studied many years ago by e.g. Todorov \cite{Todorov}, 
Olshanskii \cite{Olsh}, Molev \cite{Molev} and others. The pattern described in the above 
Proposition corresponds to the partition $p=p+0$ in the terminology of \cite{Todorov}. 
\medskip

By repeating verbatim the arguments for the compact case, one shows that the Gelfand-Tsetlin 
variables for the chain of subalgebras 
$$
\u(1)\subset \u(2)\subset \cdots \subset \u(p) \subset \u(p,1) 
\subset \cdots \subset \u(p,q-1)\subset\u(p,q) 
$$
yield a complete family of hamiltonians in involution on any admissible 
coadjoint orbit of ${\rm U}(p,q)$, which all have periodic flows. 

Now, we will show the direct derivation of the Gelfand-Tsetlin pattern, analogous to 
\cite[Proposition 6.1.3]{FM}\footnote{I thank Hermann Flaschka for explaining this to me.}
First, let ${\bf e}$ be a $2\times 2$ matrix representing an element of $\u(1,1)^*$. 
Let $\delta$ and  $\gamma$ be the eigenvalues of ${\bf e}$ with corresponding orthonormal eigenvectors
${\bf u}$ and ${\bf v}$ respectively . We assume at this moment that ${\bf u}$ is timelike and ${\bf v}$
is spacelike. They are mutually orthogonal with respect to the pseudohermitian form of 
signature $(1,1)$. Note that if ${\bf u} =(a\ b)^T$, then ${\bf v} = ({\bar b}\ {\bar a})^T$. 
Now, for a unit timelike vector ${\bf w}$ and a real number $r\in\R$, we set  
$$
L = {\bf e} + r{\bf w}\otimes{\bf w}^\dagger\ .
$$
If we decompose ${\bf w} =\alpha{\bf u} +\beta{\bf v}$, 
then we compute:
$$
\frac{\det (\lambda{\rm I}-L)}{\det(\lambda{\rm I}-{\bf e})} = 1 - r\frac{|\alpha|^2}{\lambda - \delta}
+r\frac{|\beta|^2}{\lambda-\gamma}\ .
$$
By analyzing the function in the right hand side, we see that if $\delta > \gamma$ and $r>0$, then one of its
zeroes is going to be in the interval $(-\infty, \gamma)$ and the other in $(\delta, +\infty)$.
On the contrary, if $r<0$, then the two zeroes are only possible in the interval $(\gamma, \delta)$. Note that
when $\gamma r|\alpha|^2-\delta r|\beta|^2 \le -\gamma\delta$, there indeed will be zeroes in this interval. 
This observation can directly produce the Gelfand-Tsetlin pattern, which we discuss in the next Section.

\section{Action variables}

Let us first describe the Gelfand-Tsetlin pattern for the coadjoint orbit of 
$$\Lambda = {\rm diag}(1, \underbrace{0, ..., 0}_{p-1}, 1, \underbrace{0, ...0}_{q-1}).$$ 
Let $A\in \OO_\Lambda$ be a matrix in the coadjoint orbit of $\Lambda$ and 
let $A_\ell$ be the principal left upper $\ell\times\ell$ submatrix of $A$ with real
eigenvalues $\gamma_\ell\le\delta_\ell$ complementing the $(\ell-2)$ zero eigenvalues (since the 
rank of $A_\ell$ is at most 2). Then the Gelfand-Tsetlin pattern is as follows: 
$$
\begin{array}{cccc} 
\gamma_\ell\le \gamma_{l+1}, & \delta_\ell\ge \delta_{\ell+1} & {\rm for} & p\le\ell\le n-1 \\
\gamma_\ell\ge \gamma_{\ell+1}, & \delta_\ell\le\delta_{\ell+1} & {\rm for} & 1\le\ell\le p-1 \\
{} & \gamma_\ell \le 0 & {\rm for} & \ell\le p  
\end{array}
$$
Here we think $\gamma_n=\delta_n=1$ and $\gamma_1=0$ (since we only need one 
non-trivial eigenvalue, $\delta_1$, of the $1\times 1$ matrix $A_1$).  

\medskip\noindent{\it Remark.} Note that our pattern is in agreement with the pattern (3.6) in 
\cite{Todorov}, for the decomposition $p=(p-1)+1$. \medskip

With this patter in mind, we define $M_\ell$, to be the $2\times \ell$ complex 
matrix obtained from the $n\times 2$ matrix $M=({\bf z}\ {\bf w})$ by removing the last 
$(n - \ell)$ rows. Let also $\gamma_i\le \delta_i$ be the eigenvalues of the $\ell\times\ell$
matrix $A_\ell=\eta(M_\ell)=MJ_{1,1}M^*J_{p, \ell-p}$, complementing the $(\ell-2)$ zero eigenvalues, 
which are the same as the eigenvalues of the $2\times 2$ matrix $\nu(M)=J_{1,1}M^*J_{p, \ell-p}M$.  
 
Just by analyzing the traces, one can see that $\gamma_\ell + \delta_\ell =\sum_{i=1}^\ell r_i$, 
and by repeating the arguments in (5.1) of \cite{HK}, one finds that $\delta_\ell-\gamma_\ell$
yields the length of the $\ell$-th diagonal, $d_\ell$. 

The $(n-3)$ functions on $M_{\bf r}$, namely $d_2$, ..., $d_{n-2}$, yield a completely integrable
system with periodic flows. The flow, corresponding to the hamiltonian $d_\ell$ can be visualized 
similarly to the Eucledian case, as follows. We use the action of the group ${\rm U}(1,1)$ to move the 
$(\ell+1)$-st vertex to the $t$-axis. Then we consider the $S^1$-action on the polygon, which 
revolves the vertices numbered $i+2$, ..., $n$ around the $t$-axis, while not moving all the 
other vertices (we assume, as usual, that the first vertex is placed at the origin).  

The triangle inequalities for the Minkowski space imply the following inequalities being imposed on 
the lengths $r_i$'s and $d_i$'s: 
\begin{equation}
\label{eq932}
\begin{array}{ccc}
d_\ell\ge d_{\ell-1}+r_\ell & {\rm for} & 1\le \ell \le p \\
d_\ell\ge d_{\ell+1}+r_{\ell+1} & {\rm for} & p\le \ell\le n
\end{array}
\end{equation}
Note that this is in complete agreement with the Gelfand-Tsetlin pattern described in the beginning 
of this Section. 

\section{Moduli spaces of polygons as lagrangian loci in complex quotients}

In this section we will further justify considering $M_{\bf r}$ as a reasonable geometric object. 
We will show that it is a Hausdorff topological space, and, moreover, for a generic choice
of ${\bf r}$, it has a structure of smooth manifold. In the spirit of \cite{RS}, one should always 
view quotients by real reductive groups as being homeomorphic to quotients of certain minimal loci 
by the maximal compact subgroup. In our situation, we can go a little further as in \cite{Foth_real}, 
since we have compatible involutions on groups and spaces in question at our disposal. 

Let $\tau$ stand for the complex conjugate involution of $\gl(2, \C)$ defining the real form 
$\u(1,1)$. Without any fear of confusion, we will denote by $\tau$ also the induced 
involution on the space of traceless matrices $sl(2, \C)$ as well as on the 
corresponding dual vector spaces. By using the Killing form, we identify $\su(1,1)^*$ as
the subspace of $sl(2, \C)^*$, which is the fixed point set of $\tau$. 

Let us consider the product of complex coadjoint orbits 
$\OO^\C_1\times\cdots\times\OO^\C_n$ corresponding for $1\le i \le n$ respectively to integral points 
${\rm diag}(m_i, -m_i)$ in $\ft^*$. Note that with our convention, these are fixed by $\tau$. 
The choice of $m_i$'s leads to a choice of polarization on the orbits, and therefore, 
we can consider the GIT quotient $Y=(\OO^\C_1\times\cdots\times\OO^\C_n)//{\rm SL}(2, \C)$.\
Note that the quotient map $\OO^\C_j\to \C\P^1$ by the action of the maximal unipotent subgroup $N$   
is equivariant with respect to the ${\rm SL}(2, \C)$-action and therefore
$Y$ fibers over the moduli space of eucledian polygons $P_{m_1, ..., m_n}$, which is a smooth 
projective variety, for a generic choice of ${m_i}$'s, with contractible fibers. 

The involution $\tau$ descends onto the space $Y$ and its fixed point set by \cite{RS} 
is homeomorphic to $M_{\bf r}$. 

Note that from the polygonal consideration, the isotropy subgroup of a polygon will 
be trivial if we have $\sum_{i=1}^p r_i \ne \sum_{j=p+1}^n r_j$ (in which case we can visualize 
the degenerate $n$-gon as being aligned along the $t$-axis with $p$ forward-tracks and $q$ backtracks).
Barring this situation, a polygon will represent a smooth point in the moduli space.   

\section{Final remarks}

It would be interesting to study further the topology of these moduli spaces, in particular
compute their cohomology rings. Also, one can be interested in extending these results to 
minimal elliptic orbits of more general Lie groups, the same way the Flaschka-Millson spaces \cite{FM} 
extend the moduli spaces of spatial Euclidean polygons. 

Another interesting question is to understand 
the relationship between the lattice points in the convex polyhedral set $P$ defined by \ref{eq932} 
and decompositions of the tensor product of representations of ${\rm SU}(1,1)$. In the case 
$q=1$ this question was answered by Millson in \cite{Millson}, where he showed that the 
inequalities, restricting the lengths of the diagonals, exactly match the analogue \cite{Repka} of 
the Clebsch-Gordan formula for the discrete series representations of 
${\rm SU}(1,1)\simeq {\rm SL}(2, \R)$. 

Finally, unlike in the Eucledian case, the bending torus action on $M_{\bf r}$ appears to be globally 
defined, which raises another interesting question, whether for an integral choice of ${\bf r}$, the space 
$M_{\bf r}$ has the structure of a toric variety, corresponding to the polyhedral set $P$.

\section*{Acknowledgements}

I am grateful to Hermann Flaschka for illuminating conversations, to Yi Hu for answering my question, 
and to John Millson for sending me \cite{Millson}.


\end{document}